\documentclass[12pt]{article}
\usepackage[dvips]{graphicx,epsfig}
\usepackage{amsthm}
\usepackage[centertags]{amsmath}
\usepackage{amssymb, dsfont}
\usepackage{amsfonts}
\usepackage{graphics}

\usepackage{cite}

\DeclareMathOperator{\vl}{\mathrm{vol}}

\DeclareMathOperator{\mes}{\mathrm{mes}}
\DeclareMathOperator{\sign}{\mathrm{sign}}
\def\L{\mathcal L}

\def\N{\mathbb N}
\def\Q{\mathbb Q}
\def\F{\mathbb F}
\def\B{\mathcal B}

\def\R{\mathcal R}
\def\E{\mathcal E}
\def\M{\mathcal M}

\def\la{\lambda}
\def\ka{\varkappa}

\def\vep{\varepsilon}

\def\vvvert{\vert\hspace{-1pt}\vert\hspace{-1pt}\vert}
\def\LLL{\langle\!\langle}
\def\RRR{\rangle\!\rangle}
\theoremstyle{plain}
\newtheorem{theorem}{’Ґ®аҐ¬ }
\newtheorem{lemma}{‹Ґ¬¬ }
\newtheorem{arrow}{‘«Ґ¤бвўЁҐ}
\newtoks\thehProclaim
\newtheorem*{Proclaim}{\the\thehProclaim}
\newenvironment{proclaim}[1]{\thehProclaim{#1}\begin{Proclaim}}{\end{Proclaim}}

\theoremstyle{definition}
\newtheorem{definition}{ЋЇаҐ¤Ґ«Ґ­ЁҐ}

\newcommand{\Bbbone}{\mathds 1}
\newcommand{\bp}{\begin{proof}}
\newcommand{\ep}{\end{proof}}
\newcommand{\bl}{\begin{lemma}}
\newcommand{\el}{\end{lemma}}
\newcommand{\bt}{\begin{theorem}}
\newcommand{\et}{\end{theorem}}
\newcommand{\bd}{\begin{definition}}
\newcommand{\ed}{\end{definition}}
\newcommand{\ba}{\begin{arrow}}
\newcommand{\ea}{\end{arrow}}

\begin{document}

\thispagestyle{empty}

%\rightline{POMI Preprint~-- 12/2014}
\vspace{.5in}

\vskip 1cm

\centerline{\large\bf Dyadic shift randomization in classical discrepancy theory}
%.% \footnotemark[1]{}
%${}^*$
 \vskip1cm

\centerline{\large\bf M. M. Skriganov}

  \vskip1cm

\centerline{St.~Petersburg Department of Steklov Institute of Mathematics,}

\centerline{Fontanka 27, St.~Petersburg 191023, Russia}

\vskip.1in

\centerline{E-mail: maksim88138813@mail.ru}  %skrig.@pdmi.ras.ru

 %uniform distributions, lacunar series, codiny theory

%Two mean value theorems for the $L_q$-discrepancies of point distributions

\vskip1in

\centerline{Abstract}

{\small Dyadic shifts $D\oplus T$ of point distributions $D$ in the $d$-dimensional unit cube
$U^d$ are considered as a form of randomization. Explicit formulas for the $L_q$-discrepancies of
such randomized distributions are given in the paper in terms of Rademacher functions.
Relying on the statistical independence of Rademacher functions,
Khinchin's inequalities, and other
related results, we obtain very sharp upper and lower bounds for the mean $L_q$-discrepancies,
$0<q\leq \infty$.

The upper bounds imply directly a generalization of the well known Chen's
theorem to mean discrepancies
with respect to dyadic shifts (Theorem 2.1).

From the lower bounds, it follows that for an arbitrary $N$-point
distribution $D_N$ and any exponent $0<q\leq 1$, there exist dyadic shifts
$D_N\oplus T$ such that the $L_q$-discrepancy
$\L_q[D_N\oplus T]>c_{d,q}(\log N)^{\frac{1}{2}(d-1)}$
(Theorem~2.2).

The lower bounds for the $L_{\infty}$-discrepancy are also considered in the paper. It is shown
that for an arbitrary $N$-point distribution $D_N$, there exist dyadic shifts
$D_N\oplus T$ such that $\L_{\infty}[D_N\oplus T]>c_d(\log
N)^{\frac{1}{2}d}$ (Theorem~2.3).}

\vfil


   \newpage

\centerline{\large\bf Contents}

\vspace{1cm}

1. Dyadic shifts and the mean discrepancies

2. Main results

3. Rademacher functions and related inequalities

4. Rademacher functions and explicit formulas for discrepancies

5. Explicit formulas and preliminary bounds for the mean discrepancies

6. Bounds for the error terms and some auxiliary bounds

7. Proofs of Theorem~2.1, 2.2 and 2.3

\section*{1. Dyadic shifts and the mean discrepancies}

The classical problem in discrepancy theory deals with the distribution of finite
point sets in rectangular
sub-boxes in the unit cube with sides parallel to the coordinate axes. A
detailed discussion of numerous methods
and results known in the field can be found in \cite{CST,BC,M}. We recall only
the main definitions and facts
necessary for the purposes of our paper.

Let $D$ be an arbitrary finite subset, or distribution, in the unit cube
$U^d=[0,1)^d$. The {\it local
discrepancy} $\L[D,Y]$, $Y=(y_1,\dots, y_d)\in U^d$, is defined by
\begin{equation}
\L[D,Y]=|D\cap B_Y|-|D|\vl B_Y,
\tag{1.1}
\end{equation}
where $B_Y=[0,y_1)\times \dots \times [0,y_d)$ is a rectangular box of
volume $\vl B_Y=y_1,\dots ,y_d$, and $|\cdot|$ denotes the cardinality of
a set.

The $L_q$-{\it discrepancies} are defined by
\begin{equation}
\L_q[D]=\left(\int_{U_d}|\L[D,Y]|^qdY\right)^{1/q}, \quad 0<q<\infty,
\tag{1.2}
\end{equation}
and
\begin{equation}
\L_{\infty}[D]=\sup_{Y\in U^d}|\L[D,Y]|.
\tag{1.3}
\end{equation}

We write $\N$ for the set of all positive integers, $\N_0$
for the set of all non-negative integers, $\mathbb{N}^d$ and
$\mathbb{N}^d_0$ for the product of $d$ copies of the corresponding sets.
For $s\in \N_0$, we put
$$
\Q(2^s)=\{x=m2^{-s}\in [0,1):m=0,1,\dots ,2^s-1\}
$$
and
$$
\Q^d(2^s)=\{X=(x_1,\dots, x_d)\in U_d:x_j\in \Q(2^s),j=1,\dots ,d\}.
$$
Furthermore, we put
$$
\Q(2^{\infty})=\bigcup_{s\geq 0}\Q(2^s)
\quad\mbox{and}\quad
\Q^d(2^{\infty})=\bigcup_{s\geq 0}\Q^d(2^s).
$$
The points of $\Q^d(2^{\infty})$ are called {\it dyadic rational points}.

Any $y\in [0,1)$ can be represented in the form
\begin{equation}
y=\sum_{a\geq 1}\eta_a(y)2^{-a},
\tag{1.4}
\end{equation}
where $\eta_a(y)\in \{0,1\}\simeq \F_2$, $a\in \N$. Here $\F_2$ is the
field of two elements identified with the set of residues $\{0,1\}\bmod 2$.

The dyadic expansion (1.4) is unique if we agree that for each dyadic rational point, the sum in (1.4)
contains finitely many nonzero terms. Under this convention, $\eta_a(y)=0$ for $a>s$ if $y\in \Q(2^s)$
or, in other words, for each point $y\in [0,1)$, the sequence $\{\eta_a(y):a\in \N\}$ contains
infinitely many zeros.

In a natural way, the set of dyadic rational points can be endowed with
the structure of a vector space over the finite field $\F_2$. For any two
points $x$ and $y$ in $\Q(2^{\infty})$, we define their sum $x\oplus y$ by
\begin{equation}
\eta_a(x\oplus y)=\eta_a(x)+\eta_a(y)\bmod 2,
\quad
a\in \N,
\tag{1.5}
\end{equation}
and for any two points $X=(x,\dots ,x_d)$ and $Y=(y_1,\dots ,y_d)$ in
$\Q^d(2^{\infty})$ we define
\begin{equation}
X\oplus Y=(x_1\oplus y_1,\dots ,x_d\oplus y_d).
\tag{1.6}
\end{equation}
With respect to the addition $\oplus$ defined in this way, each set
$\Q^d(2^s)$ is a vector space over the field $\F_2$, and $\dim
\Q^d(2^s)=ds$.

Note that (1.5) and (1.6) consistently define the addition $\oplus$
for all pairs of points $X$ and $Y$, whenever only one of the points, say
$Y$, belongs to $\Q^d(2^{\infty})$, while the other is an arbitrary point
$X\in U^d$.

The above shows that, for an arbitrary distribution $D$ and any point
$T\in \Q^d(2^{\infty})$, we can define the dyadic shift $D\oplus
T=\{X\oplus T:X\in D\}$ and view it as a new distribution. For each $s\in
\N$, we can consider the family $\{D\oplus T:T\in \Q^d(2^s)\}$ as a
randomization of $D$ and the corresponding discrepancies $\L_q[D\oplus T]$
as random variables.

The aim of the present paper is to study the {\it mean $L_q$-discrepancies}
\begin{equation}
\M_{s,q}[D]=\left(2^{-ds}\sum_{T\in \Q^d(2^s)}\L_q[D\oplus
T]^q\right)^{1/q}, \quad 0<q<\infty,
\tag{1.7}
\end{equation}
and
\begin{equation}
\M_{s,\infty}[D]=\max_{T\in \Q^d(2^s)}\ \L_{\infty}[D\oplus T].
\tag{1.8}
\end{equation}

Our results are given in the next section in Theorems 2.1, 2.2 and 2.3. In
Theorem~2.1, we will consider the upper bounds for $\M_{s,q}[D]$,
$0<q<\infty$, and specific distributions $D$, the so-called
$(\delta,s,d)$-nets. The lower bounds for $M_{s,q}[D]$ and arbitrary
distributions $D$ will be given in Theorems~2.2 and 2.3 for exponents
$0<q\leq 1$ and $q=\infty$ respectively.

We now recall the definition of dyadic $(\delta,s,d)$-nets.

Consider {\it elementary intervals} $\Delta^m_a\subset [0,1)$ of the
form
\begin{equation}
\Delta^m_a=[m2^{-a},(m+1)2^{-a}),
\quad
a\in \N_0\mbox{ and } m=0,1,\dots ,2^a-1,
\tag{1.9}
\end{equation}
and {\it elementary boxes} $\Delta^M_A\subset U^d$ of the form
\begin{equation}
\Delta^M_A=\Delta^{m_1}_{a_1}\times \dots \times \Delta^{m_d}_{a_d},
\quad
m_j=0,1,\dots ,2^{a_j}-1\mbox{ and } j=1,\dots,d,
\tag{1.10}
\end{equation}
where $A=(a_1,\dots, a_d)$, $M=(m_1,\dots, m_d)\in \N^d_0$.
Every such box has volume $\vl \Delta^M_A=2^{-a_1-\dots -a_d}$.

Let $0\leq \delta \leq s$ be integers. A subset $D_{2^s}\subset U^d$
consisting of $N=2^s$ points is called a {\it dyadic $(\delta,s,d)$-net of deficiency} $\delta$
if each elementary box $\Delta^M_A$ of volume $2^{\delta-s}$ contains exactly $2^{\delta}$ points of
$D_{2^s}$.

It follows from the definition that any $(\delta,s,d)$-net $D_{2^s}$ has zero discrepancy in all elementary
boxes of large volume. Precisely,
\begin{equation}
|D_{2^s}\cap \Delta^M_A|\left\{\begin{array}{ll}
=2^s\vl \Delta^M_A,&\mbox{if $\vl\ \Delta^M_A\geq 2^{\delta-s}$},\\
\leq 2^{\delta},&\mbox{if $\vl\ \Delta^M_A<2^{\delta}$}.
\end{array}\right.
\tag{1.11}
\end{equation}
Indeed, in the first case, each box $\Delta^M_A$ is a disjoint union of
elementary boxes of volume $2^{\delta-s}$, and in the second, each box
$\Delta^M_A$  is contained in an elementary box of volume $2^{\delta-s}$.

Notice also that for any $(\delta,s,d)$-net $D_{2^s}$, its shift
$D_{2^s}\oplus T$, $T\in \Q^d(2^{\infty})$, is a net with the same
parameters.

Indeed, $|(D\oplus T)\cap \Delta^M_A|=|D\cap (\Delta^M_A\oplus T)|$, $T\in
\Q^d(2^{\infty})$, and $\Delta^M_A\oplus T=\Delta^{M(T)}_A$ with an index $M(T)$.

Replacing the base 2 in the definitions (1.9) and (1.10) by an
arbitrary prime~$p$, we arrive at $(\delta,s,d)$-nets in base~$p$.
In arbitrary dimensions~$d$, the first constructions of dyadic
$(\delta,s,d)$-nets with $\delta\leq d\log d$ were given by Sobol.
Later, other constructions of nets in arbitrary base $p$ were proposed by Faure.
For details and further references, see \cite{BC,M}.

It is significant that for each base $p$, the deficiency $\delta$
increases with the growth of the dimension $d$. Furthermore,
$(0,s,d)$-nets in the base $p$ and with arbitrary large $s$ exist if and
only if $d\leq p+1$. In particular, infinite sequences of dyadic nets with
$\delta=0$ exist only in dimensions $d=1,2$ and 3.

It is known that $(\delta,s,d)$-nets $D_{2^s}$ fill the unit cube very
uniformly, and the $L_{\infty}$-discrepancies admit the bounds
\begin{equation}
\L_{\infty}[D_{2^s}]<C_d2^{\delta}s^{d-1}, \quad s\to \infty,
\tag{1.12}
\end{equation}
with a constant $C_d$ depending only on dimension $d$. Furthermore, for
arbitrary $(\delta,s,d)$-nets, the order of this bound as $s\to \infty$ cannot be improved.

We recall that for an arbitrary $N$-point distribution $D_N\subset U^d$, the bound
\begin{equation}
\L_q[D_N]>c_{d,q}(\log N)^{\frac{1}{2}(d-1)}, \quad 1<q<\infty,
\tag{1.13}
\end{equation}
holds with positive constants $c_{d,q}$ depending only on $d$ and $q$.

These classical bounds are due to Roth for $2\leq q\leq \infty$ and
Schmidt for $1<q<2$. In two dimensions, it is known that (1.13) is
also true for $q=1$, a result due to Hal\'asz.

The order of bound the (1.13) is best possible as $N\to \infty$. In the
most general form, in all dimensions $d\geq 2$ and for all exponents
$1<q<\infty$,   this fundamental fact was established by Chen. Previously,
for $1<q\leq 2$, this fact was established by Davenport, Roth and other authors.

We remark that Chen gave two different proofs of his theorem.
In the first proof \cite{C1}, averages of the $L_q$-discrepancies was considered
with respect to the usual Euclidean translations of point distributions.
The original idea of the $p$-adic shifts was introduced and exploited in
the second proof in the paper \cite{C2}.

We refer the reader to \cite{CST,BC,M} for detailed discussion of all these
questions.

The author is very grateful to William Chen for his comments and remarks on an earlier version
of this paper.

This paper was completed while the author was a visitor to the special semester
``High-dimensional Approximation'' at ICERM, Brown University.
The author thanks many participants for discussion, as well as the director and staff of ICERM
for their hospitality.

\section*{2. Main results}

Our first result concerns upper bounds for the mean $L_q$-discrepancies.

\begin{proclaim}{Theorem 2.1} Let $D_{2^s}$ be an arbitrary dyadic
$(\delta,s,d)$-net. Then, for each $0<q<\infty$, we have
\begin{equation}
\M_{s,q}[D_{2^s}]<2^{-d+\delta+1}\left(\lceil \frac{1}{2}q\rceil
(s+1)\right)^{\frac{1}{2}(d-1)}+d2^{\delta}.
\tag{2.1}
\end{equation}
In particular, there exist dyadic shifts $T\in \Q^d(2^s)$ such that
\begin{equation}
\L_q[D_{2^s}\oplus T]\leq 2^{-d+\delta+1}\left(\lceil \frac{1}{2}q\rceil
(s+1)\right)^{\frac{1}{2}(d-1)}+d2^{\delta}.
\tag{2.2}
\end{equation}
\end{proclaim}

Theorem 2.1 shows that, in all dimensions, there exist dyadic
$(\delta,s,d)$-nets which meet the lower bound (1.13).

For the first time, results of such type were established by Chen for
nets of deficiency $\delta=0$ in an arbitrary prime base $p\geq 2$.

The original Chen's approach relies on an elaborate combinatorial
analysis involving simultaneous induction on the parameters $d$, $s$, and even integers $q$.
In this approach, the assumption $\delta=0$ turns out to be essential.
As a result, for each fixed
prime base $p$, Chen's theorem could only be established in dimensions $d\leq p+1$,
and for dyadic nets only in dimensions 1, 2 and 3.
In other words, to establish Chen's theorem in dimension $d$, a prime $p\geq d-1$ needs to be chosen.

In the author's paper \cite{S1}, a new approach to the study of
the mean $L_q$-discrepancies was
proposed. In this approach, the value of the deficiency $\delta$ turns
out to be completely
irrelevant. This approach relies on the theory of lacunary function
series. In the case of
dyadic nets, these are series of Rademacher functions, which form a
lacunary subsystem of the Walsh functions.
In the case of nets in an arbitrary base $p$, these series form a lacunary subsystem of the
corresponding Chrestenson--Levy functions. The detailed description of
such functional systems can be found in \cite{GES}.

A result similar to Theorem~2.1 was established previously in \cite{S1},
see also \cite{S2}, but with worse constants in the bounds. As functions of $q$, the constants given
above in (2.1) and (2.2) are optimal in the following sense.
It can be shown that
\begin{equation}
\L_q[D_{2^s}]\leq \L_{\infty}[D_{2^s}]\leq
2^{d/\vep}\left(\L_q[D_{2^s}]+d2^{\delta+1}\right),
\tag{2.3}
\end{equation}
where $q=\vep s\to \infty$ and $\vep >0$ is an arbitrary constant, see
Lemma~6.2.

Therefore, (2.1) and (2.2) imply (1.12). Furthermore, if the
order of the constants in (2.1) and (2.2) could be improved as $q\to\infty$,
then the order of (1.12) could be also improved as $s\to\infty$
for a subsequence of $(\delta,s,d)$-nets.

Now we consider lower bounds for the mean $L_q$-discrepancies. In what
follows, $\log$ denotes the logarithm in base~2.

\begin{proclaim}{Theorem 2.2} Let $D_N\subset U^d$, $d\geq 2$, be an
arbitrary $N$-point distribution and an exponent $0<q\leq 1$ be arbitrary and
fixed. Suppose that an integer $s$ is chosen to satisfy
\begin{equation}
s\geq \log N+\frac{2d+1}{q}+\frac{1}{2}(d-1)\log (d-1)+d+1+\log d.
\tag{2.4}
\end{equation}
Then
\begin{equation}
\M_{s,q}[D_N]>\gamma_q(d) (\log N)^{\frac{1}{2}(d-1)},
\tag{2.5}
\end{equation}
where
\begin{equation}
\gamma_q(d)=2^{-(2d+1)/q-d-1}(d-1)^{-\frac{1}{2}(d-1)}.
\tag{2.6}
\end{equation}
In particular, there exist dyadic shifts $T\in \Q^d(2^s)$ such that
\begin{equation}
\L_q[D_N\oplus T]>\gamma_q(d) (\log N)^{\frac{1}{2}(d-1)}.
\tag{2.7}
\end{equation}
\end{proclaim}

Certainly, (2.5) and (2.7) hold also for $1<q<\infty$ but, in this
case, these bounds follow at once from (1.13).

In dimensions $d\geq 3$, even the exact order of the $L_1$-discrepancy is not known.
The $L_q$-discrepancies
with $0<q<1$ were never considered at all for any dimension~$d$.

Theorem~2.2 shows that, in contrast to the $L_q$-discrepancies of
individual distributions, the problem of the mean $L_q$-discrepancies can be resolved
completely for all exponents $0<q\leq 1$.

It is worth noting that Theorems~2.1 and 2.2 can be extended to the
{\it conditional mean $L_q$-discrepancies}
\begin{equation}
M_{s,q}[D,V]=\left(|V|^{-1}\sum_{T\in V}\L_q[D\oplus
T]^q\right)^{1/q}, \quad 0<q<\infty,
\tag{2.8}
\end{equation}
where $V$ is a subset in $\Q^d(2^s)$.

It turns out that the conditional means (2.8) can meet the bounds of order
(2.1) and (2.5) with very small averaging subsets $V$ of cardinality $|V|=O(s^{\omega_q(d)})$
as $s\to \infty$; here $\omega_q(d)$ is a constant independent of $s$.

Certainly, such subsets $V$ should be rather specific. Some results in
this direction were obtained in \cite{S2}, and further studies of these
intriguing questions will continue in forthcoming papers.

Our result on the mean $L_{\infty}$-discrepancy can be stated as follows.

\begin{proclaim}{Theorem 2.3} Let $D_N\subset U^d$, $d\geq 3$, be an
arbitrary $N$-point distribution. Suppose that an integer $s$ is chosen
to satisfy
\begin{equation}
s\geq \log N+\frac{1}{2}(d-2)\log (d-2)+2d+\log d.
\tag{2.9}
\end{equation}
Then
\begin{equation}
\M_{s,\infty}[D_N]>\gamma_{\infty}(d) (\log N)^{\frac{1}{2}d},
\tag{2.10}
\end{equation}
where
\begin{equation}
\gamma_{\infty}(d)=2^{-2d-1}(d-2)^{-\frac{1}{2}(d-2)}.
\tag{2.11}
\end{equation}
In particular, there exist dyadic shifts $T\in \Q^d(2^s)$ such that
\begin{equation}
\L_{\infty}[D_N\oplus T]>\gamma_{\infty}(d) (\log N)^{\frac{1}{2}d}.
\tag{2.12}
\end{equation}
\end{proclaim}

In dimensions $d\geq 3$, the exact order of the $L_{\infty}$-discrepancy
still remains an open question.

In two dimensions, the answer is known:
Schmidt's lower bound
$$
\L_{\infty}[D_N]>c\log N,
\quad
D_N\subset U^2,
$$
is best possible.

In higher dimensions, Beck's lower bound
\begin{equation}
\L_{\infty}[D_N]>c_{\vep}\log N(\log \log N)^{\frac{1}{8}-\vep},
\quad
D_N\subset U^3,
\tag{2.13}
\end{equation}
where $\vep >0$ is arbitrary small, for
three-dimensional distributions remained the only known result over many
years.
Rather recently, the stronger lower bound
\begin{equation}
\L_{\infty}[D_N]>c_d(\log N)^{\frac{1}{2}(d-1)+\eta_d},
\tag{2.14}
\end{equation}
with small constants $\eta_d\gtrsim d^{-2}$ depending only on $d$, was established in all
dimensions $d\geq 3$.
These deep results are due to Bilyk and Lacey \cite{BL1} in dimension $d=3$
and Bilyk, Lacey and Vagharshakyan \cite{BLV} in dimensions $d\geq 4$, see
also the surveys \cite{B,BL2}.

For many years, a few specialists in discrepancy theory
believes that in all dimensions $d\geq 3$, the best possible lower bound is
of the form
$$
\L_{\infty}[D_N]>c_d(\log N)^{d-1}.
$$
However, contrary to such popular belief, it was conjectured recently that the
best possible lower bound should have the form
\begin{equation}
\L_{\infty}[D_N]>c_d(\log N)^{\frac{1}{2}d}.
\tag{2.15}
\end{equation}
This latest conjecture is inspired by some very non-trivial parallels between
discrepancy theory and the theory of stochastic processes. The reader can
consult the papers \cite{BL1,BLV,B, BL2} for a more detailed discussion
of these questions.

Theorem~2.3 shows that the conjectured bound (2.15) is valid for the mean
$L_{\infty}$-discrepancy.

We will see that the mean $L_q$-discrepancies can be represented in terms
of the Rademacher series, see section~4. For such series, very sharp upper
and lower $L_q$-bounds for any $0<q<\infty$ can be given by Khinchin's
inequality. In fact, Theorems~2.1 and 2.2 are corollaries of this
inequality. At the same time, Theorem~2.3 is a corollary of a suitably
modified Khinchin's inequality, adapted to the $L_{\infty}$-norm, see
Lemma~3.2.

The lower bounds (1.13), (2.13) and (2.14) are obtained with the help of
different variations of Roth's orthogonal function method, cf. \cite{BC,B}.
It is interesting to note that, in the proofs of Theorems 2.2 and 2.3,
we will not use any auxiliary orthogonal functions. The corresponding
lower bounds (2.5) and (2.10) will be derived directly from the explicit formulas for
discrepancies given in Lemma 4.3.

\section*{3. Rademacher functions and related\\ inequalities}

In this section all necessary facts on Rademacher functions and related
topics are collected.

In the one-dimensional case, the Rademacher functions $r_a(y)$, $y\in
[0,1)$, $a\in \N$, can be defined by
\begin{equation}
r_a(y)=(-1)^{\eta_a(y)}=1-2\eta_a(y),
\tag{3.1}
\end{equation}
where $\eta_a(y)$ are the coefficients in the dyadic expansion (1.4). It
is convenient to put $r_0(y)\equiv 1$.

Notice immediately that the expansion (1.4) takes the form
\begin{equation}
y=\frac{1}{2}-\frac{1}{2}\sum_{a\geq 1}2^{-a}r_a(y).
\tag{3.2}
\end{equation}

The Rademacher functions $r_a(\cdot)$, $a\in\N$, form a sequence of
independent random variables taking the values $\pm 1$ with probability
$1/2$. This fact can be expressed by the relations
\begin{equation}
\mes\{y\in [0,1):r_{a_1}(y)=\vep_1,\dots, r_{a_l}(y)=\vep_l\}=2^{-l}
\tag{3.3}
\end{equation}
which hold for any $1\leq a_1<\dots <a_l$, $l\in \N$, and any $\vep_j=\pm
1$, $j=1,\dots ,l$, see, for example \cite{CT,PS}.

Each function $r_a(y)$, $a\in \N$, is piecewise constant on elementary
intervals $\Delta^m_a=[m2^{-a},(m+1)2^{-a})$, $m=0,1,\dots, 2^a-1$.
Therefore, the relations (3.3) are equivalent to their discrete analogs
\begin{equation}
|\{y\in \Q(2^s):r_{a_1}(y)=\vep_1,\dots, r_{a_l}(y)=\vep_l\}|=2^{s-l}
\tag{3.4}
\end{equation}
for any $1\leq a_1<\dots <a_l\leq s$, $s\in \N$, and any $\vep_j=\pm 1$,
$j=1,\dots ,l$.

The $k$-dimensional Rademacher functions $r_A(Y)$, $Y=(y_1,\dots ,y_k)\in
U^k$, $A=(a_1,\dots ,a_k)\in \N^k_0$, are defined by
\begin{equation}
r_A(Y)=\prod^d_{j=1}r_{a_j}(y_j).
\tag{3.5}
\end{equation}
In some formulas, we write $k$ for dimension, because the formulas will be
used in the subsequent text with $k=d$ and $k=d-1$.

We introduce the linear space $\R^k_s$, $s\in \mathbb{N}_0$, consisting of all
functions of the form
\begin{equation}
f(Y)=\sum_{A\in I^k_s}\lambda_Ar_A(Y),
\tag{3.6}
\end{equation}
with real coefficients $\lambda_A$. Here $I_s=\{0,1,\dots, s\}$, and
$I^k_s$ denotes the product of $k$ copies of $I_s$.

It follows from (3.4) that the set of functions
$\{r_a(\cdot):a\in I_s\}$ is linearly independent on $\Q(2^s)$, and
therefore, the set $\{r_A(\cdot): A\in I^d_s\}$ is linearly independent on
$\Q^d(2^s)$. Thus, $\dim \R^k_s=(s+1)^k$, and $\R^k_s$ is a very small
subspace in the large space $\B^k_s$ of dimension $2^{ks}$ consisting of
all real-valued functions that are piecewise constant on elementary cubes
$$
\Delta^M_s=[m_12^{-s},(m_1+1)2^{-s})\times \dots \times
[m_k2^{-s},(m_k+1)2^{-s}),
$$
where $m_j=0,1,\dots, 2^s-1$ and $ j=1,\dots, k$.

Each function $f\in \B^k_s$ is determined by its values on dyadic rational
points $\mathbb{Q}^k(2^s)$, and we have
\begin{equation}
\|f\|_q=\|f\|_{s,q}, \quad 0<q\leq \infty,
\tag{3.7}
\end{equation}
where
\begin{align*}
\|f\|_q&=\left(\int_{U^k}|f(Y)|^qdY\right)^{1/q}, \quad 0<q<\infty,
\\
\|f\|_{\infty}&=\sup_{Y\in U^k}|f(Y)|,
\\
\|f\|_{s,q}&=\left(2^{-ks}\sum_{Y\in
\Q^k(2^s)}|f(Y)|^q\right)^{1/q}, \ 0<q<\infty,
\\
\|f\|_{s,\infty}&=\max_{Y\in \Q^k(2^s)}|f(Y)|.
\end{align*}

{\it The $k$-dimensional Khinchin's inequality}: For each function
$f\in \R^k_s$ and all $0<q<\infty$,
we have
\begin{equation}
\alpha^k_qQ_2[f]\leq \|f\|_{s,q}\leq \beta^k_qQ_2[f],
\tag{3.8}
\end{equation}
{\it where}
\begin{equation}
Q_2[f]=\left(\sum_{A\in I^k_s}\lambda^2_A\right)^{1/2}.
\tag{3.9}
\end{equation}
The constants $\alpha^k_q$ and $\beta^k_q$ are independent of $f$
and $s$. They are the $k$-th powers of the constants $\alpha_q$ and $\beta_q$
respectively, with
\begin{equation}
\alpha_q\geq \left\{\begin{array}{ll}
2^{-(2-q)/q}, &\mbox{if $0<q<2$},\\
1, &\mbox{if $2\leq q<\infty$},
\end{array}\right.
\tag{3.10}
\end{equation}
and
\begin{equation}
\beta_q\leq \lceil \frac{1}{2}q\rceil^{1/2}.
\tag{3.11}
\end{equation}

In the one-dimensional case, (3.8) is a corollary of the independence of Rademacher
functions, see (3.3), (3.4). Its proof can be found in many texts on harmonic analysis and
probability theory, see, for example, \cite[Sec.~10.3, Thm.~1]{CT},
\cite{PS}, \cite[Chap.~5, Thm.~8.4]{Z}.

The extension of Khinchin's inequality to higher dimensions can be easily given by induction
on $k$; we refer the reader to \cite[Appendix~D]{S} for details.

In the subsequent text, we shall use corollaries of Khinchin's inequality
given below in Lemmas~3.1 and 3.2.

For $Y=(y_1,\dots ,y_d)\in U^d$ and $A=(a_1,\dots, a_d)\in I^d_s$, $d\geq2$, we put
\begin{equation}
\left\{
\begin{array}{ll}
Y=(\mathbf{Y},y),
&\mbox{$\mathbf{Y}=(y_1,\dots ,y_{d-1})\in U^{d-1}$ and $ y=y_d\in [0,1)$},\\
A=(\mathbf{A},a),
&\mbox{$\mathbf{A}=(a_1,\dots, a_{d-1})\in I^{d-1}_s$ and $a=a_d\in I_s$}.
\end{array}
\right.
\tag{3.12}
\end{equation}
Then any function $f\in \R^d_s$ can be written in the form
\begin{equation}
f(Y)=f(\mathbf{Y},y)=\sum_{\mathbf{A}\in
I^{d-1}_s}\Phi_{\mathbf{A}}(y)r_{\mathbf{A}}(\mathbf{Y}),
\tag{3.13}
\end{equation}
where
\begin{equation}
\Phi_{\mathbf{A}}(y)=\sum_{a\in I_s}\lambda_Ar_a(y),
\tag{3.14}
\end{equation}
as well as in the form
\begin{equation}
f(Y)=f(\mathbf{Y},y)=\sum_{a\in I_s}\varphi_a(\mathbf{Y})r_a(y),
\tag{3.15}
\end{equation}
where
\begin{equation}
\varphi_a(\mathbf{Y})=\sum_{\mathbf{A}\in
I^{d-1}_s}\lambda_Ar_{\mathbf{A}}(\mathbf{Y}).
\tag{3.16}
\end{equation}

\begin{proclaim}{Lemma 3.1} For each function $f\in \R^d_s$, we have
\begin{equation}
\|f\|_{s,q}\leq \beta^{d-1}_qQ_{\infty,2}[f], \quad 0<q<\infty,
\tag{3.17}
\end{equation}
where
\begin{equation}
Q_{\infty,2}[f]=\max_{y\in \Q(2^s)}\left(\sum_{\mathbf{A}\in
I^{d-1}_s}\Phi_{\mathbf{A}}(y)^2\right)^{1/2},
\tag{3.18}
\end{equation}
and
\begin{equation}
\|f\|_{s,q}\geq \alpha^d_qQ_2[f],
\tag{3.19}
\end{equation}
where $Q_2[f]$ is defined in {\rm(3.9)}.
\end{proclaim}

\begin{proof} Applying the right inequality in (3.8) with $k=d-1$ to
(3.13), we obtain (3.17). The bound (3.19) is just the left inequality in (3.8)
with $k=d$.
\end{proof}

Lemma~3.1 will be used in the proof of Theorems~2.1 and 2.2. For the
proof of Theorem~2.3 the following more specific result will be needed.
This result can be thought of as a modification of Khinchin's inequality
for the $L_{\infty}$-norm.

\begin{proclaim}{Lemma 3.2} For each function $f\in \R^d_s$, we have
\begin{equation}
\|f\|_{s,\infty}\geq \alpha^{d-1}_1Q_{1,2}[f],
\tag{3.20}
\end{equation}
where
\begin{align*}
Q_{1,2}[f]&=\sum_{a\in I_s}Q_2[\varphi_a],
\tag{3.21}
\\
Q_2[\varphi_a]&=\left(\sum_{\mathbf{A}\in
I^{d-1}_s}\lambda^2_A\right)^{1/2}.
\tag{3.22}
\end{align*}
\end{proclaim}

\begin{proof} First of all, we observe that the relations (3.4) imply the following
identity for each one-dimensional function $\varphi\in \R_s$.
Let
$$
\varphi(y)=\sum_{a\in I_s}\varphi_ar_a(y), \quad y\in [0,1).
$$
Then
\begin{equation}
\|\varphi \|_{s,\infty}=\sum_{a\in I_s}|\varphi_a|.
\tag{3.23}
\end{equation}

Indeed, we can assume always that $\varphi_0\geq 0$, and in view of the
relations (3.4), there exists a point $y_0\in \Q(2^s)$ such that $r_a(y_0)=\sign \varphi_a$ if
$\varphi_a\ne 0$, $a\in I_s$. Therefore
$$
\|\varphi\|_{s,\infty}\geq |\varphi(y_0)|=\sum_{a\in
I_s}|\varphi_a|.
$$
The opposite inequality is obvious, and (3.23) follows.

Applying (3.23) to (3.15), we obtain
\begin{align*}
\|f\|_{s,\infty}&=\max_{\mathbf{Y}\in
\Q^{d-1}(2^s)}\max_{y\in \Q(2^s)}|f(\mathbf{Y},y)|
\\
&=\max_{\mathbf{Y}\in \Q^{d-1}(2^s)}\sum_{a\in
I_s}|\varphi_a(\mathbf{Y})|
\\
&\geq 2^{-(d-1)s}\sum_{\mathbf{Y}\in
\Q^{d-1}(2^s)}\sum_{a\in
I_s}|\varphi_a(\mathbf{Y})|=\sum_{a\in I_s}\|\varphi_a(\cdot)\|_{s,1}
\\
&\geq \alpha^{d-1}_1\sum_{a\in
I_s}Q_2[\varphi_a]=\alpha^{d-1}_1Q_{1,2}[f],
\end{align*}
where, in the last step, we use the left inequality in (3.8) with $k=d-1$ and
$q=1$.

The proof of Lemma~3.2 is complete.
\end{proof}

\section*{4. Rademacher functions and explicit\\ formulas for discrepancies}

For an arbitrary point $y\in [0,1)$ with dyadic expansion (1.4), we denote by
\begin{equation}
y^{(s)}=\sum^s_{a=1}\eta_a(y)2^{-a}, \quad s\in \N,
\tag{4.1}
\end{equation}
its projection to $\Q(2^s)$.
For $s=0$, we put $y^{(0)}=0$, so that
\begin{equation}
y=y^{(s)}+\theta_s(y)2^{-s}, \quad s\in \N_0,
\tag{4.2}
\end{equation}
where $\theta_s(y)\in [0,1)$ for all $y\in [0,1)$.

We put
\begin{equation}
\delta^{(s)}(x,y)=\left\{\begin{array}{ll}
1, &\mbox{if $x^{(s)}=y^{(s)}$},\\
0, &\mbox{if $x^{(s)}\ne y^{(s)}$}.
\end{array}\right.
\tag{4.3}
\end{equation}

It follows immediately from (1.4) and (4.1) that the elementary intervals $\Delta^m_s$,
$m=0,1,\dots,2^s-1$, see (1.9), can be written in the form
$$
\Delta^m_s=[m2^{-s},(m+1)2^s)=\{z\in [0,1):z^{(s)}=m2^{-s}\}.
$$
Therefore
\begin{equation}
\delta^{(s)}(x,y)=\delta^{(s)}(x^{(s)}\oplus
y^{(s)})=\chi(\Delta^0_s,x^{(s)}\oplus y^{(s)})
\tag{4.4}
\end{equation}
and
\begin{equation}
\delta^{(s)}(x^{(s)}\oplus
y^{(s)})=\sum^{2^s-1}_{m=0}\chi(\Delta^m_s,x)\chi(\Delta^m_s,y).
\tag{4.5}
\end{equation}
In the sequel, we write $\chi(\E,\cdot)$ for the characteristic function of
a set $\E$. Notice that
\begin{equation}
\chi(\Delta^m_s,x)=\chi(\Delta^m_s,x^{(s)})=\chi(\Delta^m_s,x^{(a)})
\tag{4.6}
\end{equation}
for any $a\geq s$.

It follows from (4.4) and (4.5) that
$$
\delta^{(s)}(x^{(s)}\oplus y^{(s)})=\sum_{z\in
\Q(2^s)}\delta^{(s)}(x^{(s)}\oplus z)\delta^{(s)}(z\oplus y^{(s)}).
$$
Furthermore, $\delta^{(s)}(x^{(s)}\oplus y^{(s)})$ is the reproducing
kernel for the space $\B_s$; in other words,
\begin{align}
f(x)&=\sum_{y\in \Q(2^s)}\delta^{(s)}(x^{(s)}\oplus
y^{(s)})f(y)
\notag\\
&=2^s\int^1_0\delta^{(s)}(x^{(s)}\oplus y^{(s)})f(y)dy, \quad f\in
\B_s.
\tag{4.7}
\end{align}

Consider the elementary intervals
\begin{equation}
\Pi_a=\Delta^1_a=[2^{-a},2^{1-a}), \quad a\in \N.
\tag{4.8}
\end{equation}
It is convenient to put $\Pi_0=[0,1)$.

In terms of the dyadic expansion (1.4), the intervals (4.8) can be described by
\begin{equation}
\Pi_a=\{z\in [0,1):\eta_a(z)=1, \eta_i(z)=0\mbox{ for } i<a\}.
\tag{4.9}
\end{equation}
Notice that for each $s\in \N$, the set of intervals $\{\Pi_a:a>s\}$ form
a partition of the open
interval $(0,2^{-s})$.

The following result is of crucial importance in the subsequent
consideration.

\begin{proclaim}{Lemma 4.1} For each $s\in \N$, the characteristic
function $\chi([0,y),\cdot)$ of the interval $[0,y),y\in [0,1)$, has the
representation
\begin{equation}
\chi([0,y),x)=\chi^{(s)}([0,y),x)+\vep^{(s)}(x,y),
\tag{4.10}
\end{equation}
where
\begin{equation}
\chi^{(s)}([0,y),x)=\frac{1}{2}-\frac{1}{2}\sum^s_{a=1}\chi(\Pi_a,x^{(s)}\oplus
y^{(s)})r_a(y).
\tag{4.11}
\end{equation}
Furthermore, for all $x,y\in [0,1)$, we have
\begin{equation}
0\leq \chi^{(s)}([0,y),x)\leq 1
\tag{4.12}
\end{equation}
and
\begin{equation}
|\vep^{(s)}(x,y)|\leq \frac{1}{2}\delta^{(s)}(x^{(s)}\oplus y^{(s)}).
\tag{4.13}
\end{equation}
\end{proclaim}

\begin{proof} We shall check the statements of the lemma for all possible
arrangements of points $x$ and $y$.

If $x=y$, then $\chi([0,y),y)=0$, $\chi^{(s)}([0,y),y)=1/2$,
$\vep^{(s)}(y,y)=-1/2$, and the bounds (4.12) and (4.13) hold.

If $x\ne y$, we put
$$
\nu=\nu(x,y)=\min\{a\in \N:\eta_a(x)\ne \eta_a(y)\}.
$$
In view of  (4.2), we obtain
\begin{equation}
y-x=(\eta_{\nu}(y)-\eta_{\nu}(x))2^{-\nu}+(\theta_{\nu}(y)-\theta_{\nu}(x))2^{-\nu},
\tag{4.14}
\end{equation}
where $\eta_{\nu}(x)\ne \eta_{\nu}(y)$ and $0\leq
|\theta_{\nu}(y)-\theta_{\nu}(x)|<1$.
From (4.14), we conclude that

(i) $x<y$ if and only if $\eta_{\nu}(y)=1$ and $\eta_{\nu}(x)=0$;

(ii) $x>y$ if and only if $\eta_{\nu}(y)=0$ and $\eta_{\nu}(x)=1$.

Furthermore, we conclude from (4.9) that
\begin{equation}
\chi(\Pi_a,x^{(a)}\oplus y^{(a)})=\left\{\begin{array}{ll}
1, &\mbox{if $a=\nu$},\\
0, &\mbox{if $a\ne \nu$}.
\end{array}\right.
\tag{4.15}
\end{equation}

The above can be expressed by the explicit formulas
\begin{align}
\chi([0,y),x)
&
=\chi(\Pi_{\nu},x^{(\nu)}\oplus y^{(\nu)})\eta(y)
\notag\\
&
=\frac{1}{2}\chi(\Pi_{\nu},x^{(\nu)}\oplus y^{(\nu)})(1-r_{\nu}(y))
\notag\\
&
=\frac{1}{2}-\chi(\Pi_{\nu},x^{(\nu)}\oplus y^{(\nu)})r_{\nu}(y).
\tag{4.16}
\end{align}

Now, taking (4.16) and (4.15) into account, we consider the
following two possibilities:

(i) If $\nu\leq s$, then (4.10) holds with
$\vep^{(s)}(x,y)=0$, and the bounds (4.12) and (4.13) are obvious.

(ii) If $\nu>s$, then (4.10) holds with $\chi^{(s)}([0,y),x)=\frac{1}{2}$
and
$$
\vep^{(s)}(x,y)=-\frac{1}{2}\chi(\Pi_{\nu},x^{(\nu)}\oplus y^{(\nu)})r_{\nu}(y),
$$
and the bound (4.12) is obvious. The bound (4.13) holds
because $\Pi_{\nu}\subset \Delta^0_s$ and, therefore,
$$
\chi(\Pi_{\nu},x^{(\nu)}\oplus y^{(\nu)})\leq \chi(\Delta^0_s,x^{(s)}\oplus
y^{(s)})=\delta^{(s)}(x^{(s)}\oplus y^{(s)}),
$$
cf.~(4.4), (4.6).

The proof of Lemma~4.1 is complete.
\end{proof}

We emphasize that (4.16) and (4.15) imply the explicit formula
\begin{align}
\chi([0,y),x)&=\sum_{a\in \N}\chi(\Pi_a,x^{(a)}\oplus
y^{(a)})\eta_a(y)
\notag\\
&=\frac{1}{2}-\sum_{a\in \N}\chi(\Pi_a,x^{(a)}\oplus
y^{(a)})r_a(y)-\delta(x,y),
\tag{4.17}
\end{align}
where $\delta(x,y)=1$ if $x=y$ and is equal to $0$ otherwise.

Furthermore, for any $x$ and $y$ the sums in (4.17) contain at most one
nonzero term. In this sense, one can say that series in (4.17) converge for
all $x$ and $y$, while the convergence is not uniform. Lemma~4.1 shows how
we may deal with such series. Although the error terms $\vep^{(s)}$ in (4.10)
are not small, they are concentrated on small subsets.

Consider the multi-dimensional extension of the above result. For an
arbitrary point $Y=(y_1,\dots ,y_d)\in U^d$, we denote by $Y^{(s)}=(y^{(s)}_1,\dots ,y^{(s)}_d)$
its projection to $\Q^d(2^s)$, so that
$$
Y=Y^{(s)}+\Theta_s(Y)2^{-s}, \quad s\in \N_0,
$$
where
\begin{equation}
\Theta_s(Y)=(\theta_s(y_1),\dots ,\theta_s(y_d))\in U^d.
\tag{4.18}
\end{equation}

Introduce elementary boxes of the form
\begin{equation}
\Pi_A=\Pi_{a_1}\times \dots \times \Pi_{a_d}, \quad A=(a_1,\dots ,a_d)\in
\N^d_0.
\tag{4.19}
\end{equation}
Each such box has volume $\vl \Pi_A=2^{-a_1-\dots -a_d}$.

Write $\ka(A)$ for the number of nonzero elements in
$A=(a_1,\dots,a_d)\in \N^d_0$.

Multiplying (4.10) with $x=x_j$, $y=y_j$, $j=1,\dots ,d$ (recall
that $r_0(y)\equiv 1$ and $\Pi_0=[0,1))$, we obtain the following result.

\begin{proclaim}{Lemma 4.2} For each $s\in \mathbb{N}$, the characteristic function
$\chi(B_Y,X)$ of the rectangular box $B_Y=[0,y_1)\times \dots \times
[0,y_d)$, $Y\in U^d$, has the representation
\begin{equation}
\chi(B_Y,X)=\chi^{(s)}(B_Y,X)+\vep^{(s)}(X,Y),
\tag{4.20}
\end{equation}
where
\begin{equation}
\chi^{(s)}(B_y,X)=2^{-d}\sum_{A\in I^d_s}(-1)^{\ka
(A)}\chi(\Pi_A,X^{(s)})r_A(Y).
\tag{4.21}
\end{equation}
Furthermore, for all $X=(x_1,\dots ,x_d)$, $Y=(y_1,\dots ,y_d)\in U^d$, we have
\begin{equation}
0\leq \chi^{(s)}(B_Y,X)\leq 1
\tag{4.22}
\end{equation}
and
\begin{equation}
|\vep^{(s)}(X,Y)|\leq
\frac{1}{2}\sum^d_{j=1}\delta^{(s)}(x^{(s)}_j\oplus y^{(s)}_j).
\tag{4.23}
\end{equation}
\end{proclaim}

\begin{proof}  By definition
$$
\chi^{(s)}(B_y,X)=\prod^{d}_{j=1}\chi^{(s)} ([0,y_j), x_j),
$$
and (4.22) follows from (4.12).

Using (3.12), we obtain
\begin{align*}
\chi(B_Y, X) & = \chi (B _{\mathbf Y}, \mathbf X) \chi ([0,y), x)  \\
& = \chi^{(s)}(B _{\mathbf{Y}}, \mathbf X) +\varepsilon^{(s)}(\mathbf X, \mathbf Y))
(\chi^{(s)}([0,y),x)+ \vep^{(s)}(x,y)) \\
& = \chi^{(s)} (B_Y,X)+\vep^{(s)} (X,Y),
\end{align*}
where
$$
\vep^{(s)}(X,Y)=\vep^{(s)}(\mathbf X, \mathbf Y)
\chi^{(s)} ([0,y),x) +\vep^{(s)}(x,y)\chi (B_Y,X).
$$
Therefore
\begin{equation}
|\vep^{(s)} (X,Y)| \leq |\vep^{(s)}(\mathbf X,\mathbf Y)|  +
|\vep^{(s)}(x,y)|.
\tag{4.24}
\end{equation}
In the one-dimensional case, the bound (4.23) is given in  (4.13). Using (4.24),
we obtain (4.23) in all dimensions by induction on $d$.
\end{proof}

Multiplying (3.2) with $y=y_j$, $j=1,\dots,d$, we obtain
$$
y_1\dots y_d= 2^{-d} \sum_{A\in \mathbb N^d_0}
(-1)^{\ka (A)} 2^{-a_1-\dots-a_d} r_A(Y).
$$
Since $\vl B_Y=y_1\dots y_d$ and $\vl \Pi_A =2^{-a_1-\dots -a_d}$,
this can be rewritten in the form
\begin{align*}
\vl B_Y &  =2^{-d} \sum_{A\in \mathbb N^d_0} (-1)^{\ka (A)}
\vl \Pi_A \, r_A (Y)
\\
&  =\vl^{(s)} B_Y +\vep^{(s)} (Y), \quad s\in \mathbb N_0 ,
\tag{4.25}
\end{align*}
where
\begin{equation}
\vl^{(s)} B_Y =2^{-d} \sum_{A\in  I^d_s} (-1)^{\ka (A)}
\vl \Pi_A \, r_A (Y),
\tag{4.26}
\end{equation}
and $\vep^{(s)}(Y)$ satisfies the bound
\begin{equation}
|\vep^{(s)} (Y)|\leq d2^{-s-1}, \quad Y\in U^d,
\tag{4.27}
\end{equation}
easily proved by induction on $d$.

The local discrepancy (1.1) can be written in the form
\begin{equation}
\mathcal L [D,Y] =\sum_{X\in D} \mathcal L (X,Y),  \quad
\mathcal L (X,Y) = \chi (B_Y,X) -\vl B_Y.
\tag{4.28}
\end{equation}
Substituting (4.20) and (4.25)  into (4.28), we obtain
\begin{equation}
\mathcal L (X,Y)=\mathcal L^{(s)} (X,Y) +\mathcal E^{(s)} (X,Y),
\tag{4.29}
\end{equation}
where
$$
\mathcal L^{(s)}(X,Y) =2^{-2}\sum_{A\in I^d_s}
(-1)^{\ka(A)} \lambda_A (X^{(s)}\oplus Y^{(s)}) r_A (Y),
$$
$$
\lambda_A(X^{(s)}) \oplus Y^{(s)}) =\chi (\Pi_A,X^{(s)}\oplus )
-\vl \Pi_A,
$$
and
$$
\mathcal E^{(s)} (X,Y) =\vep^{(s)} (X,Y) -\vep^{(s)}(Y).
$$
In view of (4.23) and (4.27), we have
$$
|\mathcal E^{(s)} (X,Y)| \leq \frac 12
\left(\sum^{d}_{j=1} \delta^{(s)}
(x^{(s)}_j \oplus y^{(s)}_j)+d2^{-s}\right), \quad
X,Y\in U^d.
$$
For an arbitrary distribution $D\subset U^d$, we denote by
$$
D^{(s)} =\{X^{(s)} :X \in D\} , \quad s\in \mathbb N_0,
$$
its projection onto $\mathbb Q^d(2^s)$, so that $|D^{(s)}| =|D|$,
where some points of $D^{(s)}$ may coincide.

We define the {\it  micro-local  discrepancies } by
\begin{align*}
\la_A[D^{(s)}\oplus Y^{(s)}]
&
=\sum_{X\in D} \la_A (X^{(s)}\oplus Y^{(s)})
\\
&
=\sum_{X\in D}
(\chi (\Pi_A, X^{(s)} \oplus Y^{(x)}) -\vl \Pi_A)
\\
&= |(D^{(s)} \oplus Y^{(s)} ) \cap \Pi_A | - |D| \vl \Pi_A.
\tag{4.30}
\end{align*}
Substituting (4.29) into (4.28), we arrive at  the following
result summarizing the above discussion.

\begin{proclaim}{Lemma 4.3} For each   $s\in \mathbb N$, the local
discrepancy  $\mathcal L [D,Y]$ has the representation
\begin{equation}
\mathcal L[D,Y] =\mathcal L^{(s)} [D,Y] +\mathcal E^{(s)} [D,Y],
\tag{4.31}
\end{equation}
where
\begin{equation}
\mathcal L^{(s)}[D,Y] =2^{-d} \sum_{A\in I^d_s}(-1)^{\ka (A)}
\lambda_A [D^{(s)} \oplus Y^{(s)}] \, r_A (Y),
\tag{4.32}
\end{equation}
and the term  $\mathcal E^{(s)}[D,Y]$ satisfies the bound
\begin{equation}
|\mathcal E^{(s)}[D,Y]| \leq \frac 12
\left( \sum^{d}_{j=1}\delta^{(s)}_j [D^{(s)} \oplus Y^{(s)}] +d
|D| 2^{-s} \right) ,
\tag{4.33}
\end{equation}
where
\begin{equation}
\delta^{(s)}_j[D^{(s)}\oplus Y^{(s)}]  =\sum_{X\in D} \delta^{(s)}_j
(x^{(s)}_j \oplus y^{(s)}_j) .
\tag{4.34}
\end{equation}
\end{proclaim}

\section*{5. Explicit formulas and preliminary bounds\\ for the mean discrepancies}

Applying Lemma 4.3 to a shifted distribution $D\oplus T$,
$T\in \mathbb Q^d(2^s)$, we obtain
\begin{equation}
\mathcal L [D\oplus T,Y] =
\mathcal L^{(s)} [D\oplus T,Y] +\mathcal E^{(s)}
[D\oplus T,Y],
\tag{5.1}
\end{equation}
where the term $\mathcal L^{(s)} [D\oplus T,Y]$  can be written  in the
form
\begin{equation}
\mathcal L^{(s)}[D\oplus T,Y]  =\mathcal F^{(s)}[D, T\oplus Y^{(s)}, Y^{(s)}],
\tag{5.2}
\end{equation}
\begin{equation}
\mathcal F^{(s)} [D,Z,Y] =2^{-d} \sum_{A\in I^d_s}
(-1)^{\ka (A)} \lambda_A [D\oplus Z] \, r_A (Y),
\tag{5.3}
\end{equation}
and
\begin{align*}
\la_A[D\oplus Z]
&
=\sum_{X\in D}
(\chi (\Pi_A,X^{(s)}\oplus Z) - \vl \Pi_A)
\\
&
= |(D\oplus Z) \cap \Pi_A| -|D| \vl \Pi_A, \quad Z\in \mathbb Q^d (2^s).
\tag{5.4}
\end{align*}
Let $L_q(\mathbb Q^d(2^s)\times U^d)$, $0<q\leq \infty$, be the space
consisting of all functions $f(T,Y)$, $T\in \mathbb Q^d (2^s)$, $Y\in
U^d$, with $\vvvert f\vvvert_q<\infty$, where
$$
\vvvert f\vvvert_q =
\left( 2^{-ds} \sum_{T\in \mathbb Q^d(2^s)}
\int_{U^d} |f(T,Y)|^q dY\right)^{1/q}, \quad  0<q<\infty ,
$$
and
$$
\vvvert f\vvvert_{\infty} = \max_{T\in \mathbb Q^d(2^s)}
\sup_{Y\in U^d}  |f(T,Y)|.
$$
For any two functions $f_1,f_2\in L_q (\mathbb Q^d(2^s)\times U^d)$, we
have
\begin{align*}
\vvvert f_1+f_2\vvvert_q  &  \le  \vvvert f_1\vvvert_q   +  \vvvert f_2\vvvert_q , \quad
1\leq q\leq \infty  ,
\tag{5.5}
\\
\vvvert f_1+f_2\vvvert^q_q  &  \le  \vvvert f_1\vvvert^q_q   +  \vvvert f_2\vvvert^q_q , \quad
0<q\leq 1 .
\tag{5.6}
\end{align*}
For $1\leq q<\infty$, (5.5) is the standard Minkowski
inequality, while (5.6) is its modification for $0<q<1$, see
\cite[Chap.~1, (9.11), (9.13)]{Z}.

Now write
\begin{equation}
\M^{(s)}_q[D]= \vvvert \mathcal L^{(s)} [D\oplus .,.] \,\vvvert_q, \quad 0<q\leq
\infty ,
\tag{5.7}
\end{equation}
and
\begin{equation}
\mathcal E^{(s)}_q [D] = \vvvert \, \mathcal{E}^{(s)}
[D\oplus .,.] \,\vvvert_q, \quad 0<q\leq \infty .
\tag{5.8}
\end{equation}
Substituting (5.1) into (1.7)  and using (5.7),  we obtain the
upper bound
\begin{equation}
\M_{s,q} [D] \leq \M^{(s)}_q [D] +\mathcal E^{(s)}_q [{D}],\quad
1\leq q< \infty.
\tag{5.9}
\end{equation}
For $0<q\leq 1$, we can simply put
\begin{equation}
\M_{s,q}[D]\leq \M_{s,1}[D]  \leq \M^{(s)}_{1}[D] +\mathcal E^{(s)}_1 [D],
\quad 0<q\le 1.
\tag{5.10}
\end{equation}
Similarly, using (5.6), we obtain the lower bound
\begin{equation}
\M_{s,q}[D]^q\geq \M^{(s)}_q[D]^q-\mathcal E^{(s)}_q[D]^q, \quad
0<q\leq 1.
\tag{5.11}
\end{equation}

The bounds (5.9), (5.10) and (5.11) will be used in the proofs of Theorems 2.1
and 2.2.

It follows from (5.2) and (5.3) that
$\mathcal L^{(s)}[D\oplus T, Y]$, as a function of $Y\in U^d$, belongs to
the space $\mathcal B^d_s$. Hence we can use (3.7)
and write (5.7) in the form
\begin{align*}
\M^{(s)}_q[D] &  =\left(
2^{-ds} \sum_{T\in \mathbb Q^d(2^s)}
||\mathcal L^{(s)} [D\oplus T, .] ||^q_{s,q} \right) ^{1/q}
\\
&  =\left(
2^{-2ds} \sum_{T,Y\in \mathbb Q^d(2^s)}
| \mathcal L^{(s)} [D\oplus T, Y] |^q \right)^{1/q}, \quad
0<q<\infty.
\tag{5.12}
\end{align*}

The following simple observation explains  why the mean
$L_q$-discrepancies  can be expressed in terms of Rademacher series.

In the vector space of pairs
$(T,Y)\in \mathbb Q^d(2^s)\times \mathbb Q^d(2^s)\simeq \mathbb
F^{2ds}_2$, we consider the linear mapping
\begin{equation}
\tau: (T,Y)\to (T\oplus Y, Y).
\tag{5.13}
\end{equation}
Obviously,  $\tau^2=\Bbbone$, $\tau^{-1}=\tau$. Hence, $\tau$ is a one-to-one
mapping, and in the double sum in (5.12), the variables $Z=T\oplus Y$  and
$Y$ can be viewed as independent. As a result, we have
\begin{equation}
\M^{(s)}_q[D] = \left ( 2^{-ds} \sum_{Z\in \mathbb Q^d(2^s)}
\mathcal{F}^{(s)}_q[D,Z]^q\right)^{1/q}, \quad 0<q<\infty,
\tag{5.14}
\end{equation}
where
\begin{equation}
\mathcal{F}^{(s)}_q[D,Z] = \left ( 2^{-ds} \sum_{Y\in \mathbb Q^d(2^s)}
|\mathcal{F}[D,Z,Y]|^q\right)^{1/q}.
\tag{5.15}
\end{equation}

The formulas (5.14) and (5.15) will be used in
the proofs of Theorems 2.1 and 2.2.

In the case of the mean $L_{\infty}$-discrepancy, the above argument
needs to be slightly modified. First of all,  using (1.8) and
(1.3), we can write
\begin{equation}
\M_{s,\infty}[D] = \max_{T\in \mathbb Q^d(2^s)}
 \sup_{Y\in U^d} |\mathcal L [D\oplus T,Y] |
\geq \max_{T,Y\in \mathbb Q^d(2^s)}  |\mathcal L [D\oplus T,Y] |.
\tag{5.16}
\end{equation}
For $Z,Y\in \mathbb Q^d(2^s)$, we put $T=Z\oplus Y$  and
\begin{equation}
\mathcal{F}[D,Z,Y] =\mathcal L[D\oplus Z\oplus Y,Y].
\tag{5.17}
\end{equation}
With this notation, (5.1) takes the form
\begin{equation}
\mathcal{F}[D,Z,Y] = \mathcal{F}^{(s)}[D, Z,Y]+\mathcal E^{(s)}[D,Z,Y],
\tag{5.18}
\end{equation}
where $\mathcal{F}^{(s)}[D,Z,Y]$ is defined in (5.3) and
\begin{equation}
\mathcal E^{(s)} [D,Z,Y] =\mathcal E^{(s)} [D\oplus Z\oplus Y, Y].
\tag{5.19}
\end{equation}
Since $\tau$ defined in (5.13) is a one-to one mapping, we have
\begin{equation}
\max_{T,Y\in \mathbb Q^d(2^s)} |\mathcal L [D\oplus T,Y] |=
\max_{Z,Y\in \mathbb Q^d(2^s)} |\mathcal{F}[D,Z,Y] |.
\tag{5.20}
\end{equation}
This relation can be continued as follows. We have
\begin{align*}
\max_{Z,Y\in \mathbb Q^d(2^s)} |\mathcal{F}[ D,Z,Y] |
&
= \max_{Z\in \mathbb Q^d(2^s)}
\max_{Y\in \mathbb Q^d(2^s)} |\mathcal{F}[ D,Z,Y] |
\\
&
\geq 2^{-ds}  \sum_{Z\in \mathbb Q^d(2^s)}
\max_{Y\in \mathbb Q^d(2^s)} |\mathcal{F}[ D,Z,Y] |
\\
&
\geq \mathcal{F}^{(s)}_{1,\infty} [D]  -\mathcal{E}^{(s)}_{1,\infty} [D],
\tag{5.21}
\end{align*}
where
\begin{equation}
\mathcal{F}^{(s)}_{1,\infty} [D]   =2^{-ds} \, \sum_{Z\in \mathbb Q^d(2^s)} \,
\mathcal{F}^{(s)}_{\infty}  [D,Z],
\tag{5.22}
\end{equation}
\begin{equation}
\mathcal{F}^{(s)}_{\infty}[D,Z]  =  \max_{Y\in \mathbb Q^d(2^s)}
|\mathcal{F}^{(s)} [D,Z,Y] |,
\tag{5.23}
\end{equation}
and
\begin{equation}
\mathcal E^{(s)}_{1,\infty}[D] = 2^{-ds}\sum_{Z\in \mathbb Q^d(2^s)}
\, \max_{Y\in \mathbb Q^d(2^s)}  | \mathcal E^{(s)} [D,Z,Y] |.
\tag{5.24}
\end{equation}

Comparing  (5.16), (5.20) and (5.21), we obtain the lower bound
\begin{equation}
\M_{s,\infty} [D]\geq \mathcal{F}^{(s)}_{1,\infty} [D] -\mathcal
E^{(s)}_{1,\infty}[D].
\tag{5.25}
\end{equation}

This bound will be used in the proof of Theorem 2.3.

We shall call the quantities  $\M^{(s)}_q[D]$  and  $\mathcal{F}^{(s)}_{1,\infty}[D]$
the {\it principal terms}, and the quantities $\mathcal E^{(s)}_q[D]$ and
$\mathcal E^{(s)}_{1,\infty}[D]$  the {\it error terms}.

\section*{6. Bounds for the error terms and some\\ auxiliary bounds}

\begin{proclaim}{Lemma 6.1}  {\rm (i)} Let $D_{2^s}$ be an arbitrary
dyadic $(\delta,s,d)$-net. Then
\begin{equation}
\mathcal E^{(s)}_{q} [D_{2^s}] \leq  d2^{\delta},\quad 0 < q\leq \infty.
\tag{6.1}
\end{equation}

{\rm (ii)} Let $D_N\subset U^d$ be an arbitrary $N$-point distribution.
Then
\begin{equation}
\mathcal E^{(s)}_q [D_N] \leq d N2^{-s}, \quad 0<q\leq 1,
\tag{6.2}
\end{equation}
and
\begin{equation}
\mathcal E^{s}_{1,\infty}[D_N] \leq d N2^{-s}.
\tag{6.3}
\end{equation}
\end{proclaim}

\begin{proof} The functions $\delta^{(s)}_j[D^{(s)}\oplus Y^{(s)}]$,
$j=1,\dots, d$, defined in (4.34), belong to the space $\mathcal B^d_s$
and satisfy (3.7). We put
\begin{equation}
\delta^{(s)}_{j,q} [D] = \| \delta^{(s)}_j [D^{(s)}\oplus .]\|_q =
\| \delta^{(s)}_j [D^{(s)} \oplus .] \|_{s,q},
\quad 0<q\leq \infty .
\tag{6.4}
\end{equation}
Obviously,
\begin{equation}
\delta^{(s)}_{j,q}[D\oplus Z] = \delta^{(s)}_{j,q}[D], \quad Z\in
\mathbb Q^d (2^s).
\tag{6.5}
\end{equation}

Applying (4.5) to (4.34),  we obtain
\begin{equation}
\delta^{(s)}_j [D^{(s)}\oplus Z] = \sum^{2^s-1}_{m=0} N_{j,m}
\chi (\Delta^m_{s,j}, z_j),
\tag{6.7}
\end{equation}
where
$$
N_{j,m} =\sum_{X\in D} \chi (\Delta^m_s, x^{(s)}_j) =
|D \cap \Delta^m_{s,j}|,
$$
and $\Delta^m_{s,j}$ denotes the elementary box
$$
\Delta^m_{s,j}= \{X=(x_1,\dots ,x_d) \in U^d: x_j \in \Delta^m_s
\mbox{ and }
x_i\in [0,1) ,  i\ne  j  \}.
$$
Notice that $\vl \Delta^m_{s,j}=2^{-s}$. Also, for each $j=1,\dots, d$,  the boxes
$\Delta^m_{s,j}$,  $m=0,1,\dots, 2^s -1$, form a partition of the unit cube $U^d$.
Therefore
\begin{equation}
\sum^{2^s-1}_{m=0} N_{j,m} =N = |D|.
\tag{6.8}
\end{equation}

(i)  From (6.7), we obtain the bound
\begin{equation}
\delta^{(s)}_{j,q} [D] \leq \delta^{(s)}_{j,\infty} \leq
\max_m N_{j,m} , \quad 0<q\leq \infty.
\tag{6.9}
\end{equation}
Using (5.8), (4.33), and (6.5), we obtain
\begin{equation}
\mathcal E^{(s)}_q [D\oplus T] \leq \frac  12
\left(\sum^{d}_{j=1}\delta^{(s)}_{j,\infty} [D] +d|D| 2^{-s}\right),
\quad 0<q\leq \infty .
\tag{6.10}
\end{equation}
If $D_{2^s}$is an arbitrary $(\delta,s,d)$-net, then $N=2^s$ and
$N_{j,m}\leq 2^{\delta}$ for all $j$ and $m$, see (1.11).
Comparing the bounds (6.9) and (6.10) for such a net, we obtain (6.1).

(ii) From (6.7) and (6.8), we obtain the bound
\begin{equation}
\delta^{(s)}_{j,q} [D] \leq \delta^{(s)}_{j,1} [D] =
\sum^{2^s-1}_{j=1} N_{j,m}  2^{-s} =N2^{-s}, \quad 0<q\leq 1.
\tag{6.11}
\end{equation}
Using (5.8), (4.33) and (6.5), we obtain
\begin{align*}
\mathcal E^{(s)}_q [D\oplus T]
&
\leq  \mathcal E^{(s)}_1  [D\oplus T]
\\
&
\leq
\frac 12    \left( \sum^{d}_{j=1} \delta^{(s)}_{j,1} [D] +d|D| 2^{-s} \right),
\quad 0<q\leq 1.
\tag{6.12}
\end{align*}
If $D_N$ is an arbitrary $N$-point distribution, then the bounds (6.11) and
(6.12) imply (6.2).

For the function (5.19), the bound (4.33) takes the form
\begin{align*}
|\mathcal E^{(s)} [D,Z,Y]|
&
=|\mathcal E^{(s)} [D\oplus Z\oplus Y,Y]|
\\
&
\leq \frac 12 \left ( \sum^{d}_{j=1} \delta^{(s)}_j
[D^{(s)}\oplus Z]  + d|D| 2^{-s} \right),
\tag{6.13}
\end{align*}
where the right hand side is independent of $Y$.
Substituting  (6.13) into (5.24), we obtain
\begin{equation}
\mathcal E^{(s)}_{1,\infty} [D] \le \frac 12
\left( \sum^{d}_{j=1} \delta^{(s)}_{j,1} [D] +
d|D| 2^{-s} \right).
\tag{6.14}
\end{equation}
If $D_N$ is an arbitrary $N$-point distribution, then the bounds
(6.11) and (6.14) imply
(6.3).

The proof of Lemma 6.1 is complete.
\end{proof}

Next, we establish the bound (2.3) mentioned in our earlier discussion
of Theorem 2.1.

\begin{proclaim}{Lemma 6.2} For an arbitrary distribution $D\subset U^d$,
we have
\begin{equation}
\mathcal L_q[D] \le \mathcal L_{\infty} [D] \leq 2^{ds/q}
(\mathcal L_q[D] +2\mathcal E^{(s)}_{\infty}[D]), \quad
1\leq q< \infty,
\tag{6.15}
\end{equation}
where the term $\mathcal E^{(s)}_{\infty}[D]$ is defined in {\rm (5.8)}.
In particular, for an arbitrary $(\delta, s,d)$-net $D_{2^s}$ and $q=\vep
s$, $\vep >0$, the bound {\rm (6.15)} takes the form
\begin{equation}
\mathcal L_q[D_{2^s}] \leq \mathcal L_{\infty} [D_{2^s}] \leq 2^{d/\vep}
(\mathcal L_q[D_{2^s}] +d2^{\delta+1}).
\tag{6.16}
\end{equation}
\end{proclaim}

\begin{proof} It follows from (4.7) that the function
$$
\delta^{(s)}(X^{(s)} \oplus Y^{(s)}) =\prod^{d}_{j=1} \delta^{(s)}
(x^{(x)}_j \oplus y^{(s)}_j)
$$
is the reproducing kernel for the space $\mathcal B^d_s$; in other words,
\begin{equation}
f(X)=\sum_{Y\in \mathbb Q^d(2^s)} \, \delta^{(s)} (X^{(s)}
\oplus Y^{(s)} ) f(Y), \quad f\in \mathcal B^d_s.
\tag{6.17}
\end{equation}

Applying H\"older's inequality to the sum in (6.17) and taking (3.7) into
account, we obtain
\begin{align*}
\| f\|_{\infty}  & = \| f\|_{s,\infty}  \leq
\left( \sum_{Y\in \mathbb Q^d(2^s)}  |f(Y)|^q\right)^{1/q}
\\
& =  2^{ds/q} \| f\|_{s,q}= 2^{ds/q} \| f\|_q,\quad 1\leq q< \infty.
\end{align*}
In particular,
\begin{equation}
\| \mathcal L^{(s)} [D,.] \|_{\infty} \leq 2^{ds/q} \| \mathcal L^{(s)}
[D,.] \|_q,
\tag{6.18}
\end{equation}
where the $\mathcal{L}^{(s)}[D,.]$ is defined in (4.32).

On the other hand, we deduce from (4.31) and (5.8)  that
$$
\| \mathcal L^{(s)}[D,.]\|_{\infty}  \geq
\|\mathcal L [D,.] \|_{\infty} - \| \mathcal E [D,.] \|_{\infty}
\geq \mathcal L_{\infty} [D]  -\mathcal E^{(s)}_{\infty} [D]
$$
and
$$
\| \mathcal L^{(s)}[D,.]\|_q  \leq
\|\mathcal L [D,.] \|_q + \| \mathcal E [D,.] \|_\infty
\leq \mathcal L_q [D]  +\mathcal E^{(s)}_{\infty} [D].
$$

Comparing these inequalities with (6.18), we obtain
\begin{align*}
\mathcal L_{\infty} [D]  &\leq 2^{ds/q} (\mathcal L_q[D] +\mathcal
E^{(s)}_{\infty} [D])+ \mathcal E^{(s)}_{\infty} [D]
\\
& \leq  2^{ds/q} (\mathcal L_q [D] +2\mathcal E^{(s)}_{\infty} [D]).
\end{align*}
This proves the right bound in (6.15). The left bound is obvious.

If $D_{2^s}$  is a $(\delta,s,d)$-net and  $q=\vep s$, $\vep >0$, then
using the bound (6.1), we obtain (6.16).

The proof of Lemma 6.2 is complete.
\end{proof}

To conclude this section, we give one further auxiliary result that
will be used in the proofs of  Theorems 2.2 and 2.3.

Consider the subset
\begin{equation}
J^k_{\sigma}(s) =\{ \Pi_A:A\in I^k_s\mbox{ and }  \vl \Pi_A=2^{-\sigma}\},\quad
\sigma \in \mathbb N,
\tag{6.19}
\end{equation}
of the $k$-dimensional elementary boxes
$\Pi_A\subset U^k$, $k\geq 2$, see (4.19).

\begin{proclaim}{Lemma 6.3} If $s\geq \sigma$, then the subset
$J^k_{\sigma}(s)=J^k_{\sigma}$  is independent of $s$, and
\begin{equation}
|J^k_{\sigma}| \geq \left( \frac{\sigma}{k-1}\right)^{k-1}.
\tag{6.20}
\end{equation}
\end{proclaim}

\begin{proof} Since $\vl \Pi_A=2^{-a_1-\dots -a_k}$, the subset (6.19)
consists of boxes $\Pi_A$ with $A=(a_1,\dots, a_k)\in I^k_s$, where
\begin{equation}
a_1+ \dots + a_k=\sigma.
\tag{6.21}
\end{equation}
Each solution of (6.21) satisfies
$0\leq a_j\leq \min \{ \sigma,s\}$, $j=1,\dots ,k$. For $s\geq \sigma$, the set of all solutions is independent of $s$.

If $s\geq \sigma$, then for any $(a_1,\dots, a_{k-1})\in \mathbb N^{k-1}_{0}$ with
$0\leq a_j\leq \lfloor \sigma/(k-1)\rfloor$, $j=1,\dots, k-1$,
the integer
$a_k=\sigma -a_1-\dots -a_{k-1}$  satisfies $0\leq a_k\leq \sigma$. Therefore,
$A=(a_1,\dots,a_k)$  is a solution of (6.21), and
$$
|J^k_{\sigma}|\geq \left ( 1+\lfloor \sigma/(k-1)\rfloor
\right)^{k-1} \geq \left ( \frac{\sigma}{k-1} \right)^{k-1}.
$$
\end{proof}

\section*{7. Proofs of Theorems 2.1, 2.2 and 2.3}

The proof of each of Theorems 2.1, 2.2 and 2.3  consists of two steps.
First, relying on the bounds for sums of Rademacher functions given in
Lemmas 3.1 and 3.2, we establish very good bounds for the principal
terms $\M^{(s)}_q[D]$  and $\mathcal{F}^{(s)}_{1, \infty}[D]$. Next,  relying on the
upper bounds for the error terms $\mathcal E^{(s)}_q[D]$ and $\mathcal E^{(s)}_{1,\infty}$
given in Lemma 6.1, we compare the principal terms with the corresponding
mean discrepancies $\M_{s,q}[D]$.

\begin{proof}[Proof of Theorem {\rm2.1}] Let $D_{2^s}$ be a
$(\delta,s,d)$-net. We first study the quantity (5.15). Applying (3.17) to
(5.3), we have
\begin{equation}
\mathcal{F}^{(s)}_q [D_{2^s},Z] \leq \beta^{d-1}_q  \, Q_{\infty,2} [\mathcal{F}^{(s)}],
\tag{7.1}
\end{equation}
where
\begin{equation}
Q_{\infty,2}[\mathcal{F}^{(s)}] = 2^{-d} \max_{y\in \mathbb Q(2^s)}
\left( \sum_{\mathbf A \in I^{d-1}_s}  \Phi_{\mathbf A}
(Z,y)^2\right)^{1/2},
\tag{7.2}
\end{equation}
\begin{equation}
\Phi_{\mathbf A} (Z,y) =  \sum_{a\in I_s} \lambda_A
[D_{2^s} \oplus Z] r_a(y),
\tag{7.3}
\end{equation}
and the coefficients  $\lambda_A[D_{2^s}\oplus Z]$ are defined in
(5.4).

For each $Z\in \mathbb Q^d(2^s)$, the shift $D_{2^s}\oplus Z$ is also a
$(\delta,s,d)$-net, and it follows from (1.11) that
$$
\lambda_A[D_{2^s}\oplus Z] =0 \quad
\mbox{if $\vl \Pi_A \geq 2^{\delta -s}$}.
$$
The condition on volumes can be written in the form
$$
\vl \Pi_A =\vl \Pi_{\mathbf A}  \vl \Pi_a =2^{-a_1 -\dots - a_{d-1}-a}\geq
2^{\delta -s}
$$
or $a\leq s-\delta -a_1 -\dots - a_{d-1}$. Therefore, the summation in
(7.3)  is extended to
$$
s\geq a\geq l, \quad l=\max \{0, s-\delta-a_1 -\dots - a_{d-1} +1\}.
$$
The elementary boxes $\Pi_A$ are mutually  disjoint. For a given
$\mathbf A$, all the boxes $\Pi_A=\Pi_{\mathbf A}\times \Pi_a$, $s\geq a\geq
l$, are embedded in the elementary box $\Pi_{\mathbf A}\times \Delta$,
where $\Delta =\Delta^{0}_{l-1}$ if $l\geq 1$, and $\Delta =[0,1)$  if
$l=0$. In both cases,
$\vl \Pi_{\mathbf A}\times \Pi_a \leq 2^{\delta -s}$.
Hence $|(D_{2^s}\oplus Z)\cap (\Pi_{\mathbf A}\times \Delta)|\leq 2^{\delta}$
by the definition  of $(\delta,s,d)$-nets, see (1.11).

With these bounds, we now estimate the function (7.3). We have
\begin{align*}
|\Phi_{\mathbf A} (Z, y)|  & \leq \sum^{s}_{a=l}
|\lambda_A[D_{2^s}\oplus Z]|
\\
& \leq  \sum^{s}_{a=l} |(D_{2^s}\oplus  Z) \cap \Pi_A|+
|D_{2^s}| \sum^{s}_{a=l} \vl \Pi_A
\\
& \leq |(D_{2^s}\oplus Z) \cap  (\Pi_{\mathbf A}\times \Delta)| +
2^s \vl (\Pi_{\mathbf A}\times \Delta) \leq 2^{\delta+1}.
\end{align*}
Substituting this into (7.2), we obtain
$$
Q_{\infty,2} [\mathcal{F}^{(s)}] \leq 2^{-d+\delta +1}|I^{d-1}_s|^{1/2} =
2^{-d+\delta +1} (s+1)^{\frac 12 (d-1)},
$$
and therefore
$$
\mathcal{F}^{(s)}_q[D_{2^s,Z}] \leq \beta^{d-1}_q 2^{-d+\delta +1} (s+1)^{\frac12
(d-1)}.
$$
Thus, for the principal term (5.14), we have
\begin{equation}
\M^{(s)}_q [D_{2^s}]\leq 2^{-d+\delta+1} \lceil \frac 12 q\rceil^{\frac 12(d-1)}
(s+1)^{\frac 12 (d-1)},
\tag{7.4}
\end{equation}
where the bound (3.11) for the constant $\beta_q$  has also been used.

Substituting (7.4) and (6.1) into
(5.9) and  (5.10), we obtain
$$
\M_{s,q} [D_{2^s}] < 2^{-d+\delta+1}
\left( \lceil\frac 12 q\rceil (s+1) \right)^{\frac 12(d-1)} +
d2^{\delta}, \quad 0<q<\infty.
$$
The proof of Theorem 2.1 is complete.
\end{proof}

\begin{proof}[Proof of Theorem {\rm 2.2}] Let $D_N\subset U^d$, $d\geq 2$, be an $N$-point
distribution. We first study the quantity (5.15).
Applying (3.19) to (5.3), we have
\begin{equation}
\mathcal{F}^{(s)}_q[D_N,Z] \geq \alpha^d_q Q_2 [\mathcal{F}^{(s)}],
\tag{7.5}
\end{equation}
where
\begin{equation}
Q_2[\mathcal{F}^{(s)}] =2^{-d}  \left(\sum_{A\in I^d_s} \lambda_A
[D_N\oplus Z]^2 \right)^{1/2}.
\tag{7.6}
\end{equation}
The coefficients $\lambda_A[D_N\oplus Z]$ are defined in (5.4), and it is
clear that
\begin{equation}
|\lambda_A [D_N\oplus Z]| \geq
\LLL N \vl \Pi_A \RRR ,
\tag{7.7}
\end{equation}
where $\LLL t\RRR=\min \{|t-n|:n\in \mathbb Z\}$ is the distance of
$t\in \mathbb R$ from the nearest integer.
Thus
\begin{equation}
Q_2 [\mathcal{F}^{(s)}] \geq 2^{-d} \left( \sum_{A\in I^d_s} \LLL N\vl \Pi_A
\RRR^2 \right)^{1/2}.
\tag{7.8}
\end{equation}
Let $\sigma \in \mathbb N$ be chosen to satisfy
$$
2^{-2}<N2^{-\sigma} \leq 2^{-1}.
$$
Then $\LLL N \vl \Pi_A\RRR \, > 2^{-2}$ for all boxes $\Pi_A$ with $\vl
\Pi_A =2^{-\sigma}$.

Let $s\in \mathbb N$  be chosen to satisfy
\begin{equation}
s\geq \sigma =\lceil \log N\rceil +1 \geq \log N+1.
\tag{7.9}
\end{equation}
Then, using Lemma 6.3 with $k=d$, we have
$$
\sum_{A\in I^d_s}\LLL N\vl \Pi_A\RRR^2 \geq \sum_{A\in J^d_{\sigma}}
\LLL N\vl \Pi_A \RRR^2
> 2^{-4} |J^d_{\sigma} | \geq  2^{-4}
\left( \frac{\log N+1}{d-1}\right)^{d-1}.
$$
Substituting this into (7.8), we have
$$
Q_2[\mathcal{F}^{(s)}] >2^{-d-2}(d-1)^{-\frac12 (d-1)}
(\log N+1)^{\frac 12 (d-1)},
$$
and therefore
$$
\mathcal{F}^{(s)} [D_N,Z] > \alpha^d_q 2^{-d-2} (d-1)^{-\frac 12(d-1)} (\log
N+1)^{\frac 12 (d-1)}.
$$
Thus, for the principal term (5.14), we have
\begin{align*}
\M^{(s)}_q[D_N] & > c_q(d) (\log N+1)^{\frac 12(d-1)}, \quad 0<q\leq 1,
\\
c_q(d) & =2^{-2d/q -d-1} (d-1)^{-\frac 12 (d-1)},
\tag{7.10}
\end{align*}
where the bound (3.10) for the constant $\alpha_q$ has also been used.

Substituting (7.10) and (6.2) into (5.11), we have
\begin{align*}
\M_{s,q} [D_N]^q & >c_q(d)^q (\log N+1)^{\frac 12(d-1)q} -
(dN2^{-s})^q
\\
& \geq c_q(d)^q (\log N+1)^{\frac 12 (d-1)q} (1-\xi_q (s)), \quad 0<q\leq
1,
\end{align*}
where
$$
\xi_q(s) = c_q (d)^{-q} (dN2^{-s})^q.
$$

Let $s$ be chosen sufficiently large to satisfy $\xi_q(s)\leq 1/2$. To do this,
we put
$$
s\geq \log N + \frac{2d+1}{q}+ \frac 12 (d-1) \log (d-1) +d+1 +\log d,
$$
and in this case the condition (7.9) is also satisfied.

As a result,  we have
$$
\M_{s,q}[D_N]> \gamma_q(d) (\log N+1)^{\frac 12 (d-1)}, \quad 0<q\leq 1,
$$
where
$$
\gamma_q(d) =2^{-1/q} c_q (d) =2^{-(2d+1)/q-d-1}
(d-1)^{-\frac 12(d-1)}
$$
The proof of Theorem 2.2 is complete.
\end{proof}

\begin{proof}[Proof of Theorem {\rm 2.3}] Let $D_N\subset U^d$, $d\geq 3$,
be an $N$-point  distribution.
We first study the quantity (5.23).
Applying (3.20) to
(5.3), we have
\begin{equation}
\mathcal{F}^{(s)}_{\infty}[D_N,Z]\geq \alpha^{d-1}_{1} Q_{1,2} [\mathcal{F}^{(s)}],
\tag{7.11}
\end{equation}
where
\begin{equation}
Q_{1,2}[\mathcal{F}^{(s)}] =2^{-d} \sum_{a\in I_s} Q_2 [\varphi_a],
\tag{7.12}
\end{equation}
and
\begin{equation}
Q_{2}[\varphi_a] =  \left( \sum_{\mathbf A\in I^{d-1}_s}
\lambda_A [D_N\oplus Z]^2\right)^{1/2}.
\tag{7.13}
\end{equation}
Using (7.7), we deduce that
\begin{equation}
Q_2[\varphi_a] \geq \left( \sum_{\mathbf A\in I^{d-1}_s}
\LLL N \vl \Pi_A\RRR^2\right)^{1/2}.
\tag{7.14}
\end{equation}
Notice that $\vl \Pi_A=\vl \Pi_A \vl \Pi_a =\vl \Pi_A 2^{-a}$. Define
$\sigma_a\in \mathbb N$  by
$$
2^{-2}<N2^{-\sigma_a-a} \leq 2^{-1}.
$$
Then $\LLL N\vl \Pi_A\RRR >2^{-2}$ for all boxes $\Pi_{\mathbf A}$ with
$\vl \Pi_{\mathbf A}=2^{-\sigma_a}$.

It is clear that $\sigma_a=\sigma-a$, $0\leq a\leq \sigma$, where
$$
\sigma=\lceil \log N\rceil +1 \geq \log N+1.
$$
Assume now that
$$
0\leq a\leq \frac 12 \sigma \quad \text{and} \quad \sigma \geq \sigma_a
\geq \frac 12 \sigma.
$$
Let $s\in \mathbb N$  be chosen to satisfy $s\geq \sigma$. Then
\begin{equation}
s\geq \sigma =\sigma_0 > \sigma_1  >\sigma_2 > \dots.
\tag{7.15}
\end{equation}
Using Lemma 6.3 with $k=d-1$, we have
\begin{align*}
\sum_{\mathbf A\in I^{d-1}_{s}} \LLL  N\vl \Pi_A\RRR^2
&
\geq
\sum_{\mathbf A \in J^{d-1}_{\sigma_a}} \LLL N\vl \Pi_A \RRR^2
>2^{-4}|J^{d-1}_{\sigma_a}|
\\
&
\geq 2^{-4}
\left( \frac{\sigma_a}{d-2}\right)^{d-2} \geq 2^{-4}
\left( \frac{\sigma/2}{d-2}\right)^{d-2}.
\end{align*}
Hence, for the quantity  (7.13), we have
$$
Q_2[\varphi_a] >2^{-2} (d-2)^{-\frac 12(d-2)} (\sigma/2)^{\frac 12 d-1},
\quad 0\leq a\leq \sigma/2.
$$
Substituting  this into (7.12), we have
\begin{align*}
Q_{1,2} [\mathcal{F}^{(s)}] &  \geq 2^{-d} \sum_{0\leq a\leq \sigma/2} Q_2[\varphi_a]
 > 2^{-d-2}(d-2)^{-\frac 12(d-2)} (\sigma/2)^{\frac 12 d}
\\
& \geq 2^{-\frac 32 d-2} (d-2)^{-\frac 12 (d-2)} (\log N+1)^{\frac 12d},
\end{align*}
and therefore
$$
\mathcal{F}^{(s)}_{\infty} [D_N,Z] > \alpha^{d-1}_{1} 2^{-\frac 32d-2} (d-2)^{-\frac 12(d-2)}
(\log N+1)^{\frac 12 d}.
$$
Thus, for the principal term (5.22), we have
\begin{equation}
\mathcal{F}^{(s)}_{1,\infty} [D_N]  >c_{\infty} (d) (\log N+1)^{\frac 12 d},
\quad c_{\infty}(d) =2^{-2d} (d-2)^{-\frac 12 (d-2)},
\tag{7.16}
\end{equation}
where the bound (3.10) for the constant $\alpha_1$ has also been used.

Substituting (7.16) and (6.3) into
(5.25), we have
\begin{align*}
\M_{s,\infty} [D_N]  & > c_{\infty} (d) (\log N+1)^{\frac 12d}
-\frac 12 (d+1) N2^{-s}
\\
& \geq  c_{\infty}(d) (\log N+1)^{\frac 12d} (1-\xi_{\infty} (s)),
\end{align*}
where
$$
\xi_{\infty} (d) =c_{\infty} (s)^{-1} (dN2^{-s}).
$$
Let $s$ be chosen sufficiently large to satisfy $\xi_{\infty}(s)\leq \frac 12$.
To do this, we put
$$
s\geq \log N +\frac 12 (d-2) \log (d-2) +2d +\log d,
$$
and in this case the condition (7.15) is also satisfied.

As a result, we have
$$
\M_{s,\infty} [D_N] >\gamma_{\infty} (d) (\log N+1)^{\frac 12 d},
$$
where
$$
\gamma_{\infty}(d)=\frac 12 c_{\infty} (d) =2^{-2d-1}
(d-2)^{-\frac 12 (d-2)}.
$$

The proof of Theorem 2.3 is complete.
\end{proof}

\end{document}